\expandafter\def\csname ver@fixltx2e.sty\endcsname{} 
\documentclass[journal]{IEEEtran}
\usepackage{cite}
\usepackage{float}
\usepackage{graphicx} %,subcaption} %[dvips]
\usepackage{tikz}
\usepackage{pgfplots}
\pgfplotsset{compat=1.14}
\usepackage{bm}
\usepackage{color}
\usepackage{xcolor}
\usepackage{amstext,epsfig,amssymb,amsbsy,amsmath}
\usepackage{mathtools}
\usepackage[output-decimal-marker={.}]{siunitx}
\usepackage{array}
\usepackage{tabularx}
\usepackage{multirow}
\usepackage{multicol}
\usepackage{setspace}
\usepackage{enumitem}
\usepackage{dblfloatfix}
\usepackage{wrapfig}
\usepackage{url}
\usepackage{hyperref}
\usepackage{cleveref}
\usepackage{eurosym}					% euro sign

\usepackage{comment}

\expandafter\def\expandafter\UrlBreaks\expandafter{\UrlBreaks%  save the current one
  \do\a\do\b\do\c\do\d\do\e\do\f\do\g\do\h\do\i\do\j%
  \do\k\do\l\do\m\do\n\do\o\do\p\do\q\do\r\do\s\do\t%
  \do\u\do\v\do\w\do\x\do\y\do\z\do\A\do\B\do\C\do\D%
  \do\E\do\F\do\G\do\H\do\I\do\J\do\K\do\L\do\M\do\N%
  \do\O\do\P\do\Q\do\R\do\S\do\T\do\U\do\V\do\W\do\X%
  \do\Y\do\Z}

\usetikzlibrary{positioning,fit,arrows.meta,backgrounds}
\usetikzlibrary{patterns,shapes,positioning}
\usetikzlibrary{decorations.markings}
\tikzset{
    module/.style={%
        draw, rounded corners,
        minimum width=#1,
        minimum height=8mm,
        font=\sffamily
        },
    module/.default=2.1cm,
    >=LaTeX
}
\usetikzlibrary{calc}

\usepackage{ifthen}

\definecolor{gray1}{rgb}{0.23529,0.23529,0.23529}%
\definecolor{gray2}{rgb}{0.49412,0.49412,0.49412}%
\definecolor{gray3}{rgb}{0.86275,0.86275,0.86275}%
\definecolor{blue1}{rgb}{0.20392,0.30196,0.49412}%
\definecolor{blue2}{rgb}{0.72941,0.83137,0.95686}%
\definecolor{blue3}{rgb}{0.15294,0.22745,0.37255}%
\definecolor{orange1}{rgb}{0.87059,0.49020,0.00000}%
\definecolor{orange2}{rgb}{1.00000,0.60000,0.00000}%
\definecolor{orange3}{rgb}{1.00000,0.89020,0.66667}%
\definecolor{yellow1}{rgb}{1.0000,0.8353,0.3098}
\definecolor{red2}{rgb}{0.70588,0.00000,0.00000}%

\definecolor{new_green1}{rgb}{0.6667,0.8510,0.3843}%
\definecolor{new_orange1}{rgb}{0.9843,0.7490,0.2706}%
\definecolor{new_purple1}{rgb}{0.6000,0.2157,0.4039}%
\definecolor{new_blue1}{rgb}{0.0314,0.2510,0.5059}%
\definecolor{new_blue2}{rgb}{0.1686,0.5490,0.7451}%
\definecolor{new_blue3}{rgb}{0.6078,0.8510,0.8275}%

\definecolor{blue_nogrid}{rgb}{0.07451,0.62353,1.00000}%
\definecolor{orange_nolosses}{rgb}{0.85098,0.32549,0.09804}%

\pgfplotsset{
    layers/my layer set/.define layer set={
        background,
        main,
        foreground
    }{},
    set layers=my layer set,
}

\allowdisplaybreaks

\begin{document}

% with this line you avoid to have ---- if an author in the bibliography is repeated for two following references
\bstctlcite{IEEEexample:BSTcontrol}

% PAPER TITLE
\title{Loss Allocation in Joint Transmission and Distribution Peer-to-Peer Markets}%

\author{Fabio~Moret, Andrea~Tosatto,~\IEEEmembership{Student~Member,~IEEE}, Thomas~Baroche, Pierre Pinson,~\IEEEmembership{Fellow,~IEEE}
\thanks{F. Moret is with A.P. Moller - Maersk, Uptake Management, K{\o}benhavn K, Denmark (email: Fabio.moret@maersk.com).}%
\thanks{A. Tosatto and P. Pinson are with the Technical University of Denmark, Department of Electrical Engineering, Kgs. Lyngby, Denmark (emails: \{antosat,ppin\}@elektro.dtu.dk).}%
\thanks{T. Baroche is with the SATIE Laboratory located at Ecole Normale Superieure de Rennes, France (email: thomas.baroche@ens-rennes.fr)}%
\thanks{This work is partly supported by the Danish ForskEL and EUDP programmes through the Energy Collective project (grant no. 2016-1-12530) and by Innovation Fund Denmark through the multiDC project (grant no. \mbox{6154-00020B}).}%
\thanks{Submitted to "IEEE Transactions on Power Systems" on January 15, 2020 - Revised on May 6, 2020 and on August 6, 2020 - Accepted on September 13, 2020.}}%

\maketitle

%% ABSTRACT %%%%%%%%%%%%%%%%%%%%%%%%%%%%%%%%%%%%%%%%%%%%%%%%%%%%%%%%%%%%%%%%%%%%%%%%%%%%%%%%%%%%%%%%%%%%%%%%%%%%%%%%
\vspace{-1em}
\begin{abstract}
Large deployment of distribute energy resources and the increasing awareness of end-users towards their energy procurement are challenging current practices of electricity markets. A change of paradigm, from a top-down hierarchical approach to a more decentralized framework, has been recently researched, with market structures relying on multi-bilateral trades among market participants. In order to guarantee feasibility in power system operation, it is crucial to rethink the interaction with system operators and the way operational costs are shared in such decentralized markets. We propose here to include system operators, both at transmission and distribution level, as active actors of the market, accounting for power grid constraints and line losses. Moreover, to avoid market outcomes that discriminate agents for their geographical location, we analyze loss allocation policies and their impact on market outcomes and prices.
\end{abstract}

%% INDEX TERMS %%%%%%%%%%%%%%%%%%%%%%%%%%%%%%%%%%%%%%%%%%%%%%%%%%%%%%%%%%%%%%%%%%%%%%%%%%%%%%%%%%%%%%%%%%%%%%%%%%%%%
\begin{IEEEkeywords}
Peer-to-peer trading, Loss allocation, Optimal power flow, Network charges.
\end{IEEEkeywords}

% CHAPTERS %%%%%%%%%%%%%%%%%%%%%%%%%%%%%%%%%%%%%%%%%%%%%%%%%%%%%%%%%%%%%%%%%%%%%%%%%%%%%%%%%%%%%%%%%%%%%%%%%%%%%%%
\section*{Nomenclature}
\setlist[description]{font=\normalfont}
Below the list of the most important symbols in alphabetical order and grouped among sets, parameters, primal and dual variables.
\vskip 0.5em
\noindent
\textsc{Sets:}
\vskip 0.2em
\begin{description}[leftmargin=3em,style=nextline, nosep]
\item[$\mathcal{C}$] Set of nodes $c$ where TSO and DSO are connected.
\item[$\mathcal{D}$] Set of distribution nodes $r$, $s$.
\item[$\mathcal{I}$] Set of market agents $i$.
\item[$\mathcal{I}_n$] Set of market agents $i$ connected to node $n$ (or $r$).
\item[$\mathcal{L}$] Set of all transmission and distribution lines.
\item[$\mathcal{L}_n^{\rm \textsc{ac}}$] Set of AC lines $k$ connected to node $n$.
\item[$\mathcal{L}_n^{\rm \textsc{dc}}$] Set of HVDC lines $h$ connected to node $n$.
\item[$\mathcal{L}^{\rm \textsc{tso}}$] Set of lines at TSO level.
\item[$\mathcal{L}^{\rm \textsc{dso}}$] Set of lines at DSO level.
\item[$\mathcal{T}$] Set of transmission nodes $n$.
\item[$\Omega_i$] Set of trading partners $j$ for agent $i$.
\end{description}
\vskip 0.5em
\noindent
\textsc{Parameters:}
\vskip 0.2em
\begin{description}[leftmargin=3em,style=nextline, nosep]
\item[$A$] Loss allocation matrix to agents.
\item[$b^0$] Line shunt susceptance.
\item[$B$] Line susceptance matrix.
\item[$D$] Loss distribution matrix (line to nodes).
\item[$\overline{F}$] Active power limits of transmission lines.
\item[$G$] Line conductance matrix.
\item[$I^{\textsc{dc}}$] Incidence matrix of HVDC lines.
\item[$K$] Agents' capacity.
\item[$M$] Linear coefficients of line losses.
\item[$N$] Power Transfer Distribution Factors (PTDF) matrix.
\vskip 0.15em
\item[$\underline{P},\overline{P}$] Agents' max and min active power setpoints.
\vskip 0.15em
\item[$\underline{Q},\overline{Q}$] Agents' max and min reactive power setpoints.
\vskip 0.15em
\item[$Q$] Constant coefficient of line losses.
\vskip 0.15em
\item[$\overline{S}$] Apparent power limits of distribution lines.
\vskip 0.15em
\item[$\text{TF}$] Modified version of the PTDF matrix.
\vskip 0.15em
\item[$\underline{V},\overline{V}$] Maximum and minimum voltage magnitudes.
\vskip 0.15em
\item[$\underline{\Theta},\overline{\Theta}$] Maximum and minimum voltage angles.
\end{description}
\vskip 0.5em
\noindent
\textsc{Variables:}
\vskip 0.2em
\begin{description}[leftmargin=3em,style=nextline, nosep]
\item[$e^{\textsc{t}},e^{\textsc{d}}$] Exchange between TSO and DSO.
\item[$f_k$] Active power flow (line $k$).
\item[$f^{\textsc{p}}_d$] Active power flow (line $d$).
\item[$f^{\textsc{q}}_d$] Reactive power flow (line $d$).
\item[$p_i$] Active power setpoint (agent $i$).
\item[$q_i$] Reactive power setpoint (agent $i$).
\item[$t_{ij}$] Trade between agent $i$ and $j$.
\item[$v_r$] Voltage magnitude at bus $r$.
\item[$w_{l}$] Line losses (line $l$ or $d$).
\item[$w_{ij}$] Losses caused by the trade $t_{ij}$.
\item[$z_{ij}$] Net injection corresponding to the trade $t_{ij}$.
\item[$\lambda_r$] Reactive nodal price.
\item[$\pi_i$] Perceived price (agent $i$).
\item[$\theta_r$] Voltage angle at bus $r$.
\item[$\tau^{\textsc{e}}$] Price of exchange TSO-DSO.
\item[$\tau^{\textsc{l}}_{ij}$] Loss price (trade $ij$).
\vskip 0.15em
\item[$\tau^{\textsc{t}}_{ij}$] Trade price (trade $ij$).
\vskip 0.15em
\item[$\tau^{\textsc{z}}_{ij}$] Grid price (trade $ij$).
\end{description}

% \vspace{-0.5em}
\section{Introduction}\label{sec:1}

The fast development of Distributed Energy Resources (DERs) and the increasing flexibility of energy management systems for residential homes and commercial buildings have paved the way to new opportunities for energy end-users. Thanks to the recent development of Information and Communication Technology (ICT) and in line with the ongoing evolution of the economy towards sharing systems \cite{selloni}, a shift from centralized to decentralized electricity markets is expected to happen. Scenarios of how these electricity markets will take place have been proposed in \cite{parag}, ranging from peer-to-peer (P2P), community-based  and microgrid structures. Multi-lateral electricity trades were already proposed in \cite{Wu1999}, but only recently a new branch of literature has started to address how these markets would be designed and implemented in practice \cite{giotitsas, Moret2018b, Morstyn2018BilateralTrading}. However, the independent evolution of electricity markets and power systems, operated by Transmission (TSO) or Distribution (DSO) System Operators, raises concerns on the ability of these new market proposals to coexist with reliable power system management. 

Pool-based electricity markets have the advantage to gather all the information needed to find market equilibria within one platform. Consequently, TSOs can interface directly with the market operator to guarantee feasibility in power system operation. This is usually done by accounting for power flows during the market clearing, like in the case of the American locational marginal pricing system, or by redispatch actions, i.e., ex-post corrective actions, as it happens in Europe. When conceiving future electricity markets that operate in a decentralized manner, the challenges are to redefine the current interactions with TSOs and to include DSOs in the market mechanism, since most of end-users are connected to the distribution level. Several solution methods to integrate System Operators (SOs) in decentralized frameworks have been proposed, mainly including two approaches: an iterative process in which SOs accept or reject multi-lateral energy trades depending on whether network constraints are violated \cite{qin, Morstyn2018Multi-ClassPreferences, Morstyn2019DesigningFlexibility} and the usage of pre-defined network charges, e.g., based on the electrical distance between peers \cite{Baroche2019ExogenousMarketsb}. 

All the proposed solutions have two main limitations. First, they do not account for a coordination between TSOs and DSOs, a crucial challenge in modern power systems, as analyzed in \cite{Papavasiliou2018CoordinationOperations,LeCadre2019ACoordination}. Second, power losses occurring in the network are often overlooked, since their procurement costs are currently included in grid tariffs of end-users. In this regard, the current literature presents several methods for allocating losses among grid users at transmission or distribution level, as pointed out in \cite{Nikolaidis2019}. However, all these methods are only suitable for the calculation of grid tariffs, as they consist in ex-post calculations of the power flows based on market outcomes. More recently, the authors of \cite{kim2019} proposed a peer-to-peer electricity market accounting for power losses, where network charges are dynamically adapted to each trade based on its grid usage. However, no attention is given to how loss allocation policies impact fairness of market equilibria.

In this work, we argue that i) loss allocation policies directly impact market outcomes and prices, and ii) many current loss allocation policies result in an unfair distribution of losses among market participants. This happens particularly at distribution level, as shown in \cite{Nikolaidis2019, Usman2019}. Therefore, this paper proposes to extend the decentralized energy market of \cite{Baroche2019ProsumerFormulation} to consider interactions between SOs and to include market products for power losses. More specifically, the contributions of this paper are:
\begin{itemize}
    \item the extension of the unified market formulation of decentralized electricity markets proposed in \cite{Baroche2019ProsumerFormulation} to account for the interactions with multiple TSOs and DSOs.
    \item the inclusion of additional market products not only to enhance coordination among different SOs, but also to dynamically compute network charges of each energy trade.
    \item the extension of the models of system operators with loss allocation policies and the investigation of their impact on fairness of market outcomes.
\end{itemize}
Such framework allows to assess fairness of payments, which is an important aspect not considered in the existing literature, by analyzing the impact of different loss allocation policies on market participants. Moreover, the obtained market architecture considers TSOs and DSOs as independent and active agents in the market clearing, resulting in market equilibria already accounting for power grid constraints, that we call endogenous P2P electricity market.  

The remainder of the paper is structured as follows. In Section \ref{sec:2}, the optimization models of each market actors are presented and then combined into a single market clearing formulation. Section \ref{sec:3} addresses the definition of different loss allocation policies and their impact on energy trades and prices in the form of network charges. A dedicated test case and simulation results are then described in Section \ref{sec:4}. Finally, Section \ref{sec:5} gathers conclusions and perspectives regarding further works.

% \vspace{-0.5em}
\section{Market Clearing Formulation}\label{sec:2}
We first lay down the optimization problems of each market participant: namely, prosumers, system operators (both at transmission and distribution level) and the market operator. We then consider the market clearing as an equilibrium problem, where all agents negotiate their energy procurement, while maximizing their profits. Eventually, we will show that this negotiation mechanism is a convex potential game, hence there exists a Nash equilibrium that can be found by solving the equivalent optimization problem. 

Since the focus of this study is to assess the properties of the proposed market design, we find the market equilibrium by directly solving the equivalent optimization problem. Distributed optimization techniques could also be readily applied to clear the market, as largely reviewed in \cite{Molzahn2017,kargarian,hug2015}. In this way, one would simulate an actual implementation of decentralized electricity markets, where each agent individually optimizes her energy procurement and exchange information with the others to reach consensus. 

\subsection{Preliminaries and Notation}
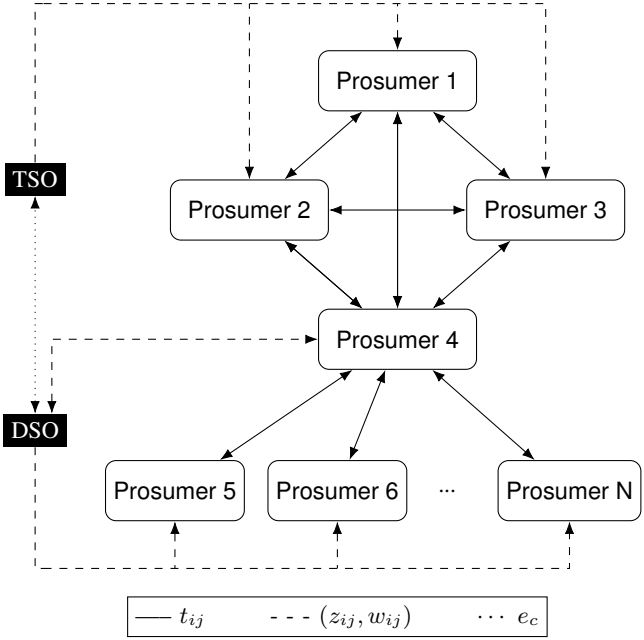
\begin{figure}[t]
    \centering
    \begin{tikzpicture}
    \small
    
    \node (I0) {};
    \node[module, above=12mm of I0] (I1) {Prosumer 1};
    \node[module, left=8mm of I0] (I3) {Prosumer 2};
    \node[module, right=8mm of I0] (I4) {Prosumer 3};
    \node[module, below=12mm of I0] (I2) {Prosumer 4};
    \node[module=1cm, below=12mm of I2, xshift=-8mm] (I5) {Prosumer 6};
    \node[module=1cm, left= 3mm of I5] (I6) {Prosumer 5};
    \node[module=1cm, right= 12mm of I5] (I7) {Prosumer N};
    \node[text width=3mm, right= 3mm of I5] (I8) {...};
    
    \node[left= 5mm of I6, yshift=8mm, fill=black] (I9) {\textcolor{white}{DSO}};
    \node[above= 29mm of I9, fill=black] (I11) {\textcolor{white}{TSO}};
    
    \node[below= 15mm of I9] (I10) {};
    \node[right= 40mm of I10] (I12) {};
    \node[above= 20mm of I11] (I13) {};
    \node[right= 40mm of I13] (I14) {};
    
    \draw[<->, dashed] (I9.45)|-(I2);
    \draw[->, dashed] (I9)|-(I10)-|(I7);
    \draw[<-, dashed] (I6) -- ($(I10)!(I6)!(I12)$);
    \draw[<-, dashed] (I5) -- ($(I10)!(I5)!(I12)$);
    
    \draw[->, dashed] (I11)|-(I13)-|(I4);
    \draw[<-, dashed] (I1) -- ($(I13)!(I1)!(I14)$);
    \draw[<-, dashed] (I3) -- ($(I13)!(I3)!(I14)$);
    
    \draw[<->, dotted] (I9)--(I11);

    \foreach \i in {1, 2, 3}
        \foreach \j in {\i, 2, 3, 4}
            \ifthenelse{\i=\j}{}{\draw[<->] (I\i)--(I\j);};
    
    \foreach \i in {5,6,7}
        \draw[<->] (I2)--(I\i);
    \draw[<->] (I2)--(I1);
    
    \node (n2) [draw,below=of I5,align=center] {----- $t_{ij}$ $\qquad$ - - - $(z_{ij}, w_{ij})$ $\qquad$ $\cdots$ $e_c$};  
\end{tikzpicture}
    \caption{Decentralized electricity market architecture where energy including grid constraints and power losses.}
    \label{fig:architecture}
    \vspace{-0.8em}
\end{figure}

The scope of this study is to integrate SOs in the negotiation process of decentralized electricity markets and to investigate properties of the consequent market outcomes. We propose a novel architecture for decentralized electricity markets, in order to account for power grid constraints and power losses. We do this by considering additional market products, on top of traded energy $t_{ij}$ between prosumer $i$ and $j$, as displayed in \figurename~\ref{fig:architecture}. These new market products include: the energy exchanged between SOs $e_c$ as well as the energy injection $z_{ij}$ and respective power losses $w_{ij}$ corresponding to $t_{ij}$. This allows for prosumers and SOs to operate as independent agents in a decentralized electricity market without the need of any information about any other market participant.

The scope of this work is to propose and analyze an innovative market architecture for decentralized electricity markets. We focus on how individual prosumers optimize their energy procurement, while interacting with independent SOs. Therefore we leave the fundamental issues of real-world implementations, power flow accuracy among all, as future work. In particular, any more accurate power flow model could be employed, as long as it maintains convexity and separability, i.e. no shared information among market agents. This justifies the need for the following assumptions. Market actors are considered rational and truthful, thus no strategic behaviour is taken into account. Any non-convexity of their optimization problems is neglected or relaxed to convex approximations, in particular the optimal power flow equations of SOs. For clarity of notation, we present the mathematical formulation of a market with one TSO and one DSO, but the extension to multiple TSOs and DSOs is straightforward: in fact, the numerical results presented in Section \ref{sec:4} consider multiple DSOs. We then restrict the market clearing to a single instance: no ramping and time coupling constraints are included in the model. Finally, we assume a perfect communication system among market actors: namely information is computed and exchanged synchronously and without any delay. This guarantees that in case the market is cleared in a decentralized fashion, the exact same market equilibrium is attained.

In each formulation, we use lower and upper case symbols respectively for variables and parameters, in bold if matrices or vectors. Dual variables, which in an equilibrium problem are decision variables of the market operator and parameters for the other agents, are expressed with Greek letters. We use $i$ as index for market participants $\mathcal{I}$, which can communicate and trade with a set of trading partners $j\in\Omega_i$. Each agent is connected to a bus of the power grid: for simplicity of notation we distinguish between TSO nodes, indexed $n \in\mathcal{T}$, and DSO buses, $r,s\in\mathcal{D}$. We identify the set of agents connected at bus $n$ (or $r$) by subscripting the set of agents $\mathcal{I}_n$. Additionally, for each bus at TSO level $n$ we split the lines connected to it in two sets, AC $k\in\mathcal{L}_n^{\textsc{ac}}$ and HVDC $h\in\mathcal{L}_n^{\textsc{dc}}$, while $\mathcal{L}^{\textsc{tso}}$ is the set of all lines at TSO level. Lines at DSO level are grouped in the set of the pairs of receiving and sending nodes $d=(r,s)\in\mathcal{L}^{\textsc{dso}}$. All the lines are then grouped in the set $\mathcal{L}$. Finally, we define a set of nodes $c\in\mathcal{C}$, where DSOs and TSOs are connected.

\subsection{Prosumers}
The objective of each prosumer $i$ is to minimize their energy procurement costs and can be modelled as an optimization problem,
\begin{subequations}\label{eq:2_1}
    \begin{alignat}{3}
	    \begin{split}
        \underset{\Gamma_{\textsc{p}}}{\text{min}} \enspace 
        & \rlap{$f_i(p_i, q_i) - \sum_{j \in \Omega_i} (\tau_{ij}^{\textsc{t}} + \tau_{ij}^{\textsc{z}}) t_{ij} - (\tau_{ij}^{\textsc{z}} + \tau_{ij}^{\textsc{l}}) w_{ij}$} \\
        & - \lambda_{r|i\in\mathcal{I}_r} \, q_i
        \end{split} \label{2_1:obj} \\
        \text{s.t.} \enspace & \sum_{j \in \Omega_i} (t_{ij} + w_{ij}) = p_i \qquad \; \qquad \qquad \qquad && : \pi_i \label{2_1:balance}\\
        & \underline{P}_i \leq p_i \leq \overline{P}_i && : \underline{\gamma}_i^{\textsc{p}},\overline{\gamma}_i^{\textsc{p}}  \label{2_1:limitP}\\
        & \underline{Q}_i \leq q_i \leq \overline{Q}_i && : \underline{\gamma}_i^{\textsc{q}},\overline{\gamma}_i^{\textsc{q}}  \label{2_1:limitQ}\\
        & w_{ij} \ge 0, t_{ij} \in \mathbb{R} \label{2_1:feas}
    \end{alignat}
\end{subequations}
where $\Gamma_{\textsc{p}} = \{ p_i, q_i, \bm{t}_i, \bm{w}_i \}$ is the set of variables, with $p_i$ the net active power of agent $i$ (positive if generated and negative if consumed) and $q_i$ the reactive power. The vector $\bm{t}_i$ is the collection of $t_{ij}$, i.e., the energy traded between agent $i$ and each of the $j$ trading partners in $\Omega_i$ at price $\tau_{ij}^{\textsc{t}}$. As in \cite{Baroche2019ProsumerFormulation}, by changing the communication topology among market participants, one could simulate different market structures, from pool-based to community-based and peer-to-peer markets.

Additionally, each trade impacts the power grid both in term of congestions, reflected by the price charge $\tau_{ij}^{\textsc{z}}$, and by creating power losses. These losses $w_{ij}$, caused by the trade between agents $i$ and $j$, need to be procured at price $\tau_{ij}^{\textsc{z}} + \tau_{ij}^{\textsc{l}}$, i.e. grid plus loss price, with the respective system operator. The sum of the trades balances out with the net generation and the losses as in \eqref{2_1:balance} at a perceived price $\pi_i$. $\overline{P}_i$ and $\underline{P}_i$ are respectively the maximum and minimum active power set-points while $\underline{\gamma}_i$ and $\overline{\gamma}_i$ are the dual variables associated with these limits. Respectively, \eqref{2_1:limitQ} bounds reactive power, whose grid costs are enforced by means of the reactive nodal price $\lambda_r$. The objective function \eqref{2_1:obj} is to minimize procurement costs, as a sum of generation costs (and consumption utility), energy and loss trades and grid costs.

\subsection{Transmission System Operator}
We formulate the optimization problem of the Transmission System Operator (TSO) as
\begin{subequations}\label{eq:2_2}
    \begin{alignat}{4}
        \begin{split}
            \underset{\Gamma_{\textsc{t}}}{\text{min}} \enspace & \rlap{ $\sum\limits_{c \in \mathcal{C}} \tau_c^{\textsc{e}} e_c^{\textsc{t}} + \sum\limits_{n \in \mathcal{T}} \sum\limits_{i \in \mathcal{I}_n} \tau_{ij}^{\textsc{z}} z_{ij}$} \\
            & \rlap{$+ \sum\limits_{l \in \mathcal{L}^{\textsc{tso}}} \sum\limits_{i | \mathcal{I}_n \in \mathcal{T}} \tau_{ij}^{\textsc{l}} A_{(i,j),l} \; w_l$}
        \end{split} \label{2_2:obj} \\
        \text{s.t.} \enspace & \rlap{$f_{k} = \sum\limits_{n \in \mathcal{T}} N_{kn} \Big( \sum\limits_{i \in \mathcal{I}_n} \; z_{ij} - \sum\limits_{h \in \mathcal{L}_n^{\textsc{dc}}} I^{\textsc{dc}}_{nh} \; f_h $} \nonumber \\
        & \quad \; - \sum\limits_{c \in \mathcal{C}_n} e_c^{\textsc{t}} - \sum\limits_{l \in \mathcal{L}} D_{nl} \; w_l  \Big) \quad &&  \forall k\in\mathcal{L}^{ac} \, && :\, \varphi_k \label{2_2:flow} \\ %\llap{ $ $}
        & -\overline{F}_l \leq f_l \leq \overline{F}_l && \forall l\in\mathcal{L}^{\textsc{tso}} && :\, \underline{\mu}_l,\overline{\mu}_l \label{2_2:limitflow} \\
        & w_l = M_l^{\textsc{t}} \; |f_l| + Q_l^{\textsc{t}} \; \; \; && \forall l\in\mathcal{L}^{\textsc{tso}} && :\, \phi_l \label{2_2:nodalLT}  \\
        & \boldsymbol{z}, \boldsymbol{e}^{\textsc{t}} \in \mathbb{R}, \; \boldsymbol{w} \ge 0 \label{2_2:feas}
    \end{alignat}
\end{subequations}
where $\Gamma_{\textsc{t}} = \{ \boldsymbol{f}, \boldsymbol{w}, \boldsymbol{z}, \boldsymbol{e}^{\textsc{t}} \}$ is the set of decision variables of the TSO, with $\boldsymbol{f}$ the vector of line flows both in AC and DC, $\boldsymbol{w}$ the vector of line losses and $\boldsymbol{e}^{\textsc{t}}$ the vector of power exchanges from TSO to DSOs. The vector $\boldsymbol{z}$ gathers all the injections corresponding to each trade, including the associated losses. The objective function \eqref{2_2:obj} represents the cost of the energy exchanged with all DSOs at price $\tau_{c}^{\textsc{e}}$, the sum of the congestion rent over AC and HVDC lines and the cost of losses. The subset $\mathcal{C}_n$ is used to identify all the DSOs connected to node $n$. On the one hand, the congestion rent is expressed as the sum over all nodes of all injections $z_{ij}$ times the respective grid price  $\tau_{ij}^{\textsc{z}}$ . On the other hand, the cost of losses is calculated as the sum over the losses of each line $l$ at price $\tau_{ij}^{\textsc{l}}$. 

Power flows over AC lines are derived using the \textit{Power Transfer Distribution Factors (PTDF) matrix}, represented by the elements $N_{kn}$ of $\bm{N}$ in \eqref{2_2:flow}, and the nodal injection of each trade, with $\boldsymbol{\varphi}$ being the respective dual variable. Nodal injections are calculated as the sum of the injections corresponding to each energy trade, HVDC line flows, energy exchange with DSOs and line losses. Line losses are transformed into nodal losses by allocating half of the line losses to the respective receiving and sending node, $D_{nl} = 0.5$ if node $n$ is connected line $l$ and zero otherwise. 

We define AC flows by means of the PTDF matrix and not through voltage angle differences, in order to avoid multiple dual solutions. In fact, in the B-theta formulation, power flows are calculated as the difference of voltage angles: the corresponding KKTs represent an undetermined system of equations with infinite dual variable solutions. The flows over HVDC lines are included in the nodal injections. Indeed, thanks to the flexibility of HVDC lines, those flows can be arbitrarily defined by the TSO and are modeled as positive or negative injections. We include HVDC lines in our TSO model since their ability to control power flows makes such components even more relevant in view of a future electricity market where power losses are included as a market product. Constraint \eqref{2_2:limitflow} defines the line limits, with $\overline{F}_l$ being the capacities of the lines, and $\underline{\mu}_l$ and $\overline{\mu}_l$ the respective dual variables, in other words the shadow prices of congestion.

In order to define the losses allocated to each trade, we first model line losses $w_l$ as a linear function of the line flows $f_l$, constraint \eqref{2_2:nodalLT}. Loss functions are included in the form of two inequality constraints as in \cite{Tosatto2019}; in this way it is possible to properly calculate losses without considering the direction of the flow (that would require the introduction of binary variables). The coefficients of the linear approximation are computed using the least squares method as in \cite{Tosatto2019}, where the authors propose an extension to piece-wise linear approximation that is added to this model for the test case analysis. Line losses are then allocated to market participants, connected at TSO $i | \mathcal{I}_n \in \mathcal{T}$, by means of allocation policies, modelled as coefficients of the matrix $\boldsymbol{A}$. By doing this, we propose loss allocation policies that directly affect the negotiation mechanism of decentralized electricity markets, instead of being allocated ex-post. We will deeper investigate different policies and their impact on market equilibria in the following sections. 

\subsection{Distribution System Operator}
A distribution system operator participates to the negotiation mechanism with objectives similar to a TSO. Given the lower voltage level, we extend the power flow equations to include reactive power, resulting in:
\\
\vspace{-0.2em}
\begin{subequations}\label{eq:2_4}
    \begin{alignat}{4}
        \begin{split}
            \underset{\Gamma_{\textsc{d}}}{\text{min}} \enspace & \rlap{ $\sum\limits_{r \in \mathcal{D}} \sum\limits_{s \in \mathcal{D}_r} \lambda_r \; f_{(r,s)}^{\textsc{q}} - \tau_c^{\textsc{e}} e_c^{\textsc{d}} + \sum\limits_{r \in \mathcal{D}} \sum\limits_{i \in \mathcal{I}_r} \tau_{ij}^{\textsc{z}} z_{ij}$} \\
            & \rlap{$+ \sum\limits_{d \in \mathcal{L}^{\textsc{dso}}} \sum\limits_{i | \mathcal{I}_r \in \mathcal{D}} \tau_{ij}^{\textsc{l}} A_{(i,j),d} \; w_d$}
        \end{split} \label{2_4:obj} \\
        \text{s.t.} \enspace &\rlap{$f_d^{\textsc{p}} = B_{rs} (\vartheta_r - \vartheta_s) - G_{rs} (v_r - v_s)$} \nonumber \\
        & \qquad \qquad  && \llap{ $\forall d:(r,s) $}\in\mathcal{L}^{\textsc{dso}} && :\, \varphi^{\textsc{p}}_{d} \label{2_4:flowP} \\
        &\rlap{$f_d^{\textsc{q}} = B^*_{rs} (v_r - v_s) + G_{rs} (\vartheta_r - \vartheta_s) - b^0_{rs}$} \nonumber \\
        & && \llap{ $\forall d:(r,s) $}\in\mathcal{L}^{\textsc{dso}}&& :\, \varphi^{\textsc{q}}_d \label{2_4:flowQ} \\
        & (f_d^{\textsc{p}})^2 + (f_d^{\textsc{q}})^2 \leq (\overline{S}_d)^2 && \llap{ $\forall d$}\in\mathcal{L}^{\textsc{dso}} && :\, \eta_d^{\textsc{ac}} \label{2_4:limitflow} \\
        & \rlap{$\sum\limits_{i,j \in \mathcal{I}_r} z_{ij} = \sum\limits_{d:(r,s)\in\mathcal{I}_r} f_d^{\textsc{p}} + \sum\limits_{d\in\mathcal{L}^{\textsc{dso}}} D_{rd} w_d + \sum\limits_{c\in\mathcal{C}_r} e_c^{\textsc{d}}$} \nonumber \\
        & &&  \forall r \in\mathcal{D} &&  :\, \eta_r \label{2_4:balance} \\
        & e_c^{\textsc{d}} = - \sum\limits_{r \in \mathcal{D}} \sum\limits_{i \in \mathcal{I}_r} z_{ij} + \sum\limits_{d \in \mathcal{L}^{\textsc{dso}}} w_d && && :\, \eta_c^{\textsc{e}} \label{2_4:flowDT} \\
        & \vartheta_{\text{ref}} = 0 \\
        & \underline{\Theta}_r \le \vartheta_r \le \overline{\Theta}_r && \forall r\in\mathcal{D} && :\, \underline{\eta}_r^{\theta}, \overline{\eta}_r^{\theta} \label{2_4:limitT}\\
        & \underline{V}_r \le v_r \le \overline{V}_r && \forall r\in\mathcal{D} && :\, \underline{\eta}^{\textsc{v}}_{\,r}, \overline{\eta}^{\textsc{v}}_{\,r} \label{2_4:limitV}\\
        & w_d = M_d^{\textsc{d}} \; |f_d^{\textsc{p}}| + Q_d^{\textsc{d}} \qquad \qquad \quad && \llap{ $\forall d$}\in\mathcal{L}^{\textsc{dso}} && :\, \phi_d \label{2_4:nodalLD} \\
        & \boldsymbol{f}^{\textsc{p}}, \boldsymbol{f}^{\textsc{q}}, \boldsymbol{z}, e^{\textsc{t}} \in \mathbb{R}, \; \boldsymbol{w} \ge 0 \label{2_4:feas}
    \end{alignat}
\end{subequations}
\\
\noindent
with $\Gamma_{\textsc{d}} = \{ \boldsymbol{f}^{\textsc{p}}, \boldsymbol{f}^{\textsc{q}}, \boldsymbol{z}, \boldsymbol{w}, e^{\textsc{d}}, \boldsymbol{\vartheta}, \boldsymbol{v} \}$ is the set of decision variables of the DSO, with $\boldsymbol{f}^{\textsc{p}}$ and $\boldsymbol{f}^{\textsc{q}}$ respectively the line active and reactive flows, $\boldsymbol{z}$ the injections corresponding to agent trades and losses, $\boldsymbol{w}$ the line losses, $e^{\textsc{d}}$ the energy flow from DSO to TSO and $\boldsymbol{\vartheta}, \boldsymbol{v}$ respectively the voltage angles and magnitudes. As already presented for the TSO model, the objective function \eqref{2_4:obj} represents the sum of the congestion rent (both for active and reactive flows), the costs of trades with the TSO and the cost of losses. The losses are traded between market participants, at distribution level $i | \mathcal{I}_r \in \mathcal{D}$, and the DSO, who enforces loss allocation policies similarly to the TSO problem, summarized with the matrix $\boldsymbol{A}$, further discussed in the following section. 

The active and reactive power flows are derived through constraints \eqref{2_4:flowP} and \eqref{2_4:flowQ} using a linear approximation of an AC power flow \cite{Yang2019}, where $b^0_{rs}$ is the shunt susceptance, $G_{rs}$ the conductance and $B_{rs}$ the susceptance of line between $d$ nodes $r$ and $s$ ($B^*_{rs} = B_{rs} +2b^0_{rs}$), with $\varphi^{\textsc{p}}_d$ and $\varphi^{\textsc{q}}_d$ the associated dual variables. Line capacities are enforced by the conic constraint \eqref{2_4:limitflow}, while nodal balances are guaranteed by \eqref{2_4:balance}.

Similarly to the formulation of the TSO problem, power losses in the distribution grid are calculated for each line through constraint \eqref{2_4:nodalLD}, by means of a linear approximation. The energy flow between TSO and DSO is calculated as the sum of all the trades originated at DSO level and the losses on all DSO lines, as in constraint \eqref{2_4:flowDT}. The case of one connection node is presented above only for simplicity of indexing and notation, however several connection points can be readily considered. Finally, constraints \eqref{2_4:limitV} and \eqref{2_4:limitT} ensure voltage standard compliance at each node $r$, with $\underline{V}_r$, $\overline{V}_r$, $\underline{\Theta}_r$ and $\overline{\Theta}_r$ being the maximum deviation of magnitudes and angles, and $\underline{\eta}^{\textsc{v}}_{\,r}$, $\overline{\eta}^{\textsc{v}}_{\,r}$, $\underline{\eta}^{\theta}_{\,r}$ and $\overline{\eta}^{\theta}_{\,r}$ the corresponding dual variables.

\subsection{Market Operator}
The market constraints ensure trade reciprocity of energy and losses, as well the allocation policies for losses and power balances at TSO and DSO level. The respective optimization problem becomes
\begin{subequations}\label{eq:2_5}
    \begin{alignat}{2}
    \underset{\Gamma_{\textsc{m}}}{\text{min}} \enspace & \sum\limits_{i\in\mathcal{I}} \sum\limits_{j\in\Omega_i} \Big[ -(t_{ij} + t_{ji}) \tau_{ij}^{\textsc{t}} - (t_{ij} + w_{ij} - z_{ij}) \tau_{ij}^{\textsc{z}} \nonumber \\
    & - \left(w_{ij} - \sum\limits_{l \in \mathcal{L}} A_{(i,j),l} \; w_l \; \tau_{ij}^{\textsc{l}} \right) \tau_{ij}^{\textsc{l}} \Big] +  \sum\limits_{c\in\mathcal{C}} (e_c^{\textsc{t}} - e_c^{\textsc{d}}) \tau_c^{\textsc{e}} \nonumber \\
    & + \sum\limits_{r \in \mathcal{D}} \sum\limits_{s \in \mathcal{D}_r} ( f_{(r,s)}^{\textsc{q}} - \sum\limits_{i \in \mathcal{I}_r} q_i ) \lambda_r \label{2_5:obj} \\
    \text{s.t.} \enspace & \boldsymbol{\tau}^{\textsc{t}}, \boldsymbol{\tau}^{\textsc{z}}, \boldsymbol{\tau}^{\textsc{e}}, \boldsymbol{\tau}^{\textsc{l}}, \boldsymbol{\lambda} \in \mathbb{R}
    \end{alignat}
\end{subequations}
where $\Gamma_{\textsc{m}} = \{ \boldsymbol{\tau}^{\textsc{t}}, \boldsymbol{\tau}^{\textsc{z}}, \boldsymbol{\tau}^{\textsc{e}}, \boldsymbol{\tau}^{\textsc{l}}, \boldsymbol{\lambda} \}$ is the set of decision variables of the market operator, indicating respectively the prices of the energy trades among market participants, of grid usage for each injection caused by the energy trades, of energy exchange between TSO and DSOs, of losses and the nodal prices for reactive power. The objective function of the market operator \eqref{2_5:obj} is to enforce shared constraints by setting all the mentioned prices.

\subsection{Equivalent Optimization Problem} \label{sec:eq_opt}
The optimization problems of all market players are convex problems, as long as the cost or utility function $f_i$ of each prosumer $i$ is convex. Note that a usual formulation of a DC power flow via the PTDF matrix involves employing the sum of agent injections $p_i$ in constraint \eqref{2_2:flow}. This, however, results in a generalized Nash equilibrium, since prosumer actions limit the feasible space of the TSO. By introducing the additional variables $\boldsymbol{z}$, no variable of a market player appears in the constraints of other agents, refining the generalized Nash equilibrium to a Nash equilibrium. It is straightforward to show that the Jacobian of the game map is symmetric, hence the game map is integrable (Theorem 1.3.1 in \cite{facchinei2007finite}). Therefore, a Nash equilibrium exists and we can write the following equivalent optimization problem,
\begin{subequations}\label{eq:2_6}
    \begin{alignat}{4}
    \underset{\Gamma_{\textsc{p}}, \Gamma_{\textsc{t}}, \Gamma_{\textsc{d}}}{\text{min}} \enspace & \sum\limits_{i \in \mathcal{I}} f_i(p_i, q_i) \label{2_6:obj} \\
    \text{s.t.} \enspace \enspace & \eqref{2_1:balance} - \eqref{2_1:feas} \label{2_6:pros}\\
    & \eqref{2_2:flow} - \eqref{2_2:feas} \label{2_6:TSO}\\
    & \eqref{2_4:flowP} - \eqref{2_4:feas} \label{2_6:DSO}\\
    & t_{ij} = {\color{blue}-} t_{ji} && \mathllap{\forall i\in\mathcal{I}, \;} \forall j\in\Omega_i \; && :\; \tau_{ij}^{\textsc{t}} \label{2_6:trade_rec} \\
    & t_{ij} + w_{ij} = z_{ij} && \mathllap{\forall i\in\mathcal{I}, \;} \forall j\in\Omega_i \; && :\; \tau_{ij}^{\textsc{z}} \label{2_6:inj_bal} \\
    & w_{ij} = \sum\limits_{l \in \mathcal{L}} A_{(i,j),l} \; w_l && \mathllap{\forall i\in\mathcal{I}, \;} \forall j\in\Omega_i \; && :\; \tau_{ij}^{\textsc{l}} \label{2_6:loss_bal} \\
    & e_c^{\textsc{t}} = e_c^{\textsc{d}}  && \forall c\in\mathcal{C} && :\; \tau_c^{\textsc{e}} \label{2_6:TSO_DSO_bal} \\
    & \sum\limits_{s\in\mathcal{D}_r} f_{rs}^{\textsc{q}} - \sum\limits_{i \in \mathcal{I}_r} q_i = 0 \qquad \quad \; \; && \forall r\in\mathcal{D} && :\; \lambda_r^{\textsc{q}} \label{2_6:react_bal}
    \end{alignat}
\end{subequations}
where the objective function and constraint \eqref{2_6:pros} include the optimization problems of all prosumers, and constraints \eqref{2_6:TSO} and \eqref{2_6:DSO} describe the problems of respectively the TSO and the DSO. All these agent based problems are linked by constraints \eqref{2_6:trade_rec}-\eqref{2_6:react_bal}, corresponding to the objective functions of all market agents, where such constraints were relaxed in their Lagrangian form.

Note that \eqref{eq:2_6} can be solved by means of decomposition techniques, such as the Alternating Direction Method of Multipliers (ADMM). This decentralized negotiation mechanism would read as an extension of the ADMM algorithm in \cite{Baroche2019ProsumerFormulation}, with \eqref{2_6:trade_rec} - \eqref{2_6:react_bal} the respective complicating constraints. The presence of SOC constraints in a subproblem, as in the case of \eqref{2_4:limitflow}, has been shown to not undermine the properties of ADMM since it preserves convexity and compactness of the feasible space \cite{ma2016consensus}.
\section{Loss Allocation Policies}\label{sec:3}
When including losses in decentralized electricity markets, fairness of agent payments may be undermined. To verify this, we analyze the price formation mechanism and the consequent effect of different loss allocation policies.

\subsection{Price Formation}

Trade prices, for both energy and losses, can be extracted as dual variables of \eqref{2_6:trade_rec} - \eqref{2_6:react_bal} or can be analytically derived by computing the first-order conditions of \eqref{eq:2_6}. We hereby presents only the most meaningful among all these equations to understand the impact of grid constraints and losses on trade prices. These are
\begin{subequations} \label{eq:KKT}
    \begin{alignat}{2}
        p_i :\; & \frac{\partial f_i}{\partial p_i} - \pi_i + \overline{\gamma}_i^{\textsc{p}} - \underline{\gamma}_i^{\textsc{p}} = 0 \label{eq:KKT_p} \\
        t_{ij} :\; & \pi_i - \tau_{ij}^{\textsc{t}} - \tau_{ij}^{\textsc{z}} = 0 \label{eq:KKT_t} \\
        w_{ij} :\; & \pi_i - \tau_{ij}^{\textsc{z}} - \tau_{ij}^{\textsc{l}} = 0 \label{eq:KKT_w} \\
        z_{ij} :\; & \tau_{ij}^{\textsc{z}} - \sum\limits_{k\in\mathcal{L}} N_{kn|n\in\mathcal{I}_n\cup\mathcal{C}} \; \varphi_k = 0 \label{eq:KKT_z} \\
        w_l :\; & \sum\limits_{i\in\mathcal{I}} A_{(i,j),l} \; \tau_i^{\textsc{l}} + \phi_l + \sum\limits_{k\in\mathcal{L}} N_{kn} \; D_{nl} \; \varphi_k = 0 \label{eq:KKT_w_l}
    \end{alignat}
\end{subequations}

From \eqref{eq:KKT_p}, we can identify $\pi_i$ as the price perceived by agent $i$, driving her dispatch decision, limited by asset capacities. From \eqref{eq:KKT_t} and \eqref{eq:KKT_w}, it follows that the perceived price equals each trade price $\tau_{ij}^{\textsc{t}}$, or equivalently the loss price $\tau_{ij}^{\textsc{l}}$, adjusted by $\tau_{ij}^{\textsc{z}}$, including the impact of all grid constraints on the market equilibrium as in \eqref{eq:KKT_z} (it is intuitive that a price difference happens only if some of these constraints are binding). This implies that energy and losses are traded at the same grid-aware price. Therefore, the allocation policies, summarized by $A_{(i,j),l}$, do not only impact the definition of the loss price of each trade, as in \eqref{eq:KKT_w_l}, but also the energy price of each trade. Hence, loss allocation policies should be carefully designed in order to grant fairness of market participation. Since the concept of fairness can differ among system operators, we propose two policies, that can be independently employed by different SOs or even combined to obtain a milder strategy.

\subsection{Loss Allocation Policies}
We summarize the allocation policy coefficients in a matrix $\boldsymbol{A}$ with $2\,T$ rows, with $T$ the number of trades (we define a coefficient for the trade $(i,j)$ and one for the trade $(j,i)$) and  with $L$ columns, with $L$ the number of lines.
\begin{comment}We use as reference case a socialization policy, where each trade is allocated the same share of the total system losses. $$A_{(i,j),l}^{\textsc{ref}} = \frac{1}{2T} \quad \forall i \in \mathcal{I}, \forall l \in \mathcal{L} $$
\end{comment}
We propose two different allocation policies: one where each System Operator (SO) socializes the losses in its lines ($\boldsymbol{A}^{\textsc{soc}}$) and one reflecting grid usage of individual trades ($\boldsymbol{A}^{\textsc{ind}}$). We define the socialization policy coefficients as
\begin{equation} \label{eq:soc_pol}
A_{(i,j),l}^{\textsc{soc}} = 
\begin{cases}
    \frac{1}{2T^{\textsc{so}}} \quad \quad & \text{if } i\in\mathcal{I}^{\textsc{so}}, l \in \mathcal{L}^{\textsc{so}} \vspace{0.5em} \\
    0 \quad \quad & \text{otherwise} 
    \end{cases}
\end{equation}
with $T^{\textsc{so}}$ the number of trades of the SO, that agent $i$ belongs to. Each agent is equally allocated the losses occurring on the lines of the respective SO. This policy erases geographical discrimination of market participants, but at the same time trades are not affected by their individual impact on the grid. 

The second proposed policy attains the opposite effect. By allocating to each trade the losses produced by the flows it generates, this policy discriminates agents based on their geographical location. The coefficients of this policy are defined as
\begin{equation} \label{eq:ind_pol}
A_{(i,j),l}^{\textsc{ind}} = \frac{|\text{TF}_{ln|i\in\mathcal{I}_n} - \text{TF}_{lm|j\in\mathcal{I}_m}|}{\sum\limits_{i,j\in\mathcal{I}} | \text{TF}_{ln|i\in\mathcal{I}_n} - \text{TF}_{lm|j\in\mathcal{I}_m} |}
\end{equation}
where $\text{TF}_{ln|i\in\mathcal{I}_n}$ is a modified version of the PTDF matrix for line $l$ and node $n$ where agent $i$ is connected. The PTDF matrix is modified to include DSO nodes and lines. This is done by defining the PTDF matrix for each block of nodes and lines of each SO. Then each node at DSO has the same values for the TSO lines as the node of connection between TSO and DSO. To compute the allocation policy coefficients for trade $(i,j)$, we employ the difference of the modified PTDF at the respective nodes in absolute value, to avoid negative values while accounting for grid usage at the same time. Finally, the values are normalized to guarantee that all the coefficients for each line sum to 1.

Ideally, this individual allocation policy should account for the amount of energy traded, to allocate more losses to larger trades and vice-versa.  However, this would make appear bilinear terms in the optimization problem of the SOs, which is highly non-practical. Therefore, we propose to account for the size of the market agents, by adjusting the allocation policy coefficients by the agent capacities ($K_i$). The policy in \eqref{eq:ind_pol} then becomes
\begin{equation}\label{eq:cap_pol}
A_{(i,j),l}^{\textsc{ind}} = \frac{|\text{TF}_{ln|i\in\mathcal{I}_n} - \text{TF}_{lm|j\in\mathcal{I}_m}|}{\sum\limits_{i,j\in\mathcal{I}} | \text{TF}_{ln|i\in\mathcal{I}_n} - \text{TF}_{lm|j\in\mathcal{I}_m} |} \; \frac{K_i}{\sum\limits_{i\in\mathcal{I}} K_i}
\end{equation}
where
\begin{equation}
    K_i = \text{max}\big\{|\underline{P}_i|,|\overline{P}_i|\big\}
\end{equation}

Arguments for fairness and unfairness of both policies exist. On the one hand, a socialization policy equally distributes total loss costs across market participants without discrimination on geographical location. At the same time, agents are not incentivized to trade locally, as they are charged the same no matter the distance with their trading partners. On the other hand, an individual policy makes sure that each agent pays for the losses she causes in the system, incentivizing local trading. However, this policy may discriminate, for instance, one agent at the end of a DSO feeder that gets allocated the losses of all the lines to reach the feeder, compared to a market participant close to the connection with the TSO. Therefore, we propose a linear combination of the two strategies,
\begin{equation}
A_{(i,j),l} = \chi \; A_{(i,j),l}^{\textsc{soc}} + (1-\chi) \; A_{(i,j),l}^{\textsc{ind}}
\end{equation}
with $\chi$ a socialization factor, to find a trade-off between the advantages and caveats of both policies. One should notice that different SOs can have different allocation policies. In fact in the following section, we investigate a test case with both TSO and DSOs over different combination of allocation policies to show how they impact market outcomes and payments.

\subsection{Illustrative Example} \label{sec:illustrative_example}

\begin{figure}[!b]
	\centering
	\includegraphics[trim= 0.1cm 0.5cm 0.1cm 1.1cm,clip,width=0.375\textwidth]{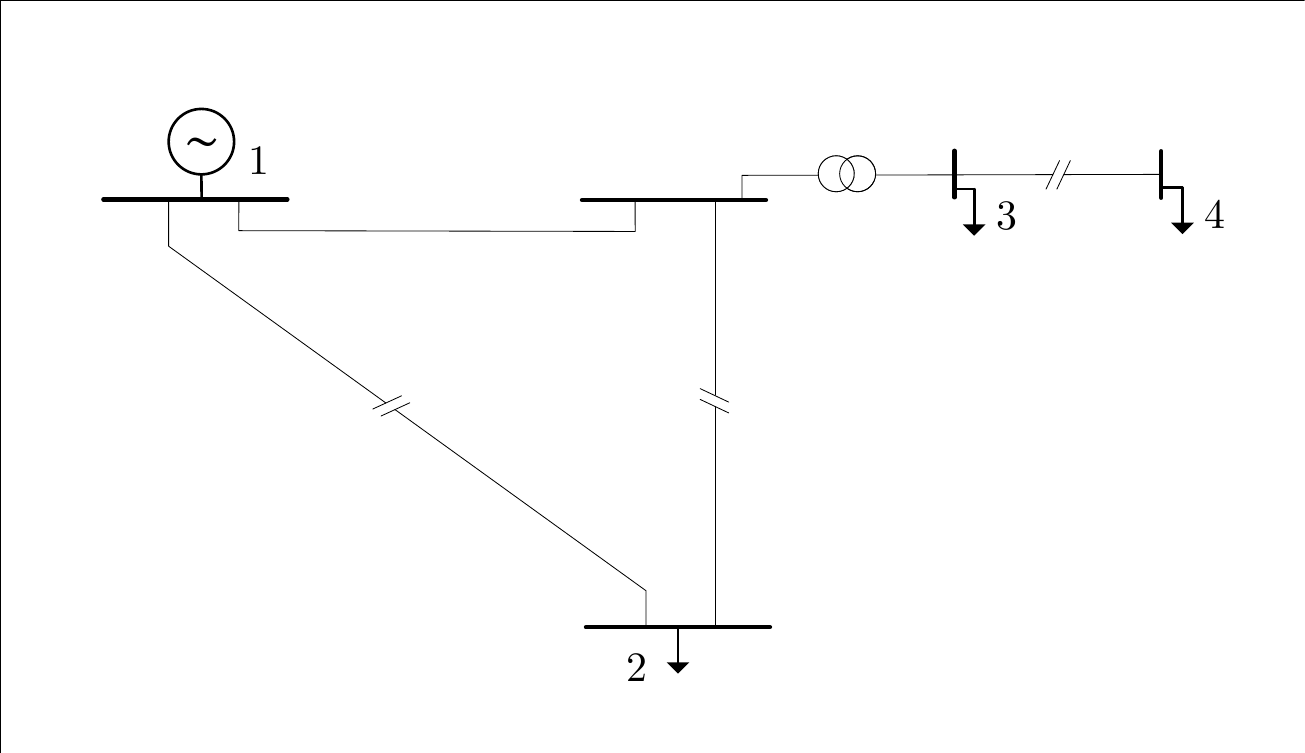}
	\vspace{-1em}
	\caption{Illustrative 5-bus test case, with a TSO and a DSO.}%
	\label{fig:4_small_case}%
	%\vspace{-1em}
\end{figure}

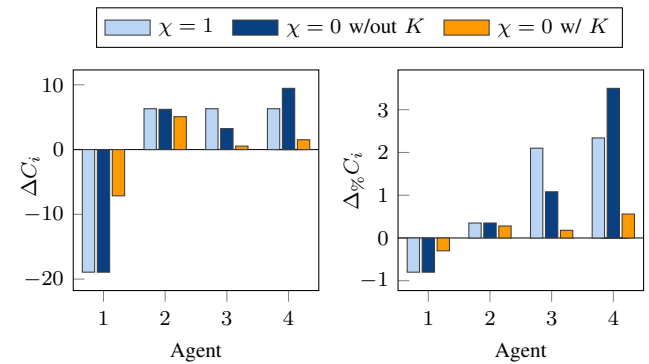
\begin{figure}[!b]
    \begin{tikzpicture}
            \begin{axis}[%
                %compat=1.3,
                width=0.18\textwidth,
                height=0.12\textheight,
                at={(1.20in,0.806in)},
                scale only axis,
                %bar shift auto,
                ybar=1pt,
                log origin=infty,
                xmin=0.51,
                xmax=4.49,
                xtick={1,2,3,4},
                xtick pos=left,
                ytick pos=left,
                xlabel style={font=\color{white!15!black}},
                xlabel={Agent},
                ylabel style={font=\color{white!15!black}, at={(-0.1,0.52)}},  % , rotate=-90
                ylabel={$\Delta C_i$},
                %ylabel near ticks,
                label style={font=\footnotesize},
                every tick label/.append style={font=\footnotesize},
                axis background/.style={fill=white},
                legend columns=3,
                legend style={at={(-1.30,1.1)}, anchor=south west, legend cell align=left, align=left, draw=black, font=\footnotesize},
                every axis legend/.append style={column sep=0.1em},
                ]
                
                \addplot[black] table[x index = 0, y index = 1, forget plot] {
                0.5 0
                4.49 0
                };
                
                \draw[ultra thin] (axis cs:\pgfkeysvalueof{/pgfplots/xmin},0) -- (axis cs:\pgfkeysvalueof{/pgfplots/xmax},0);
                
                \addplot[ybar=0pt, bar width=0.2, fill=blue2, draw=gray1, area legend] table[row sep=crcr] {%
                1 -18.93\\ 
                2 6.31\\
                3 6.31\\ 
                4 6.31\\
                };
                \addplot[ybar=0pt, bar width=0.2, fill=new_blue1, draw=gray1, area legend] table[row sep=crcr] {%
                1 -18.93\\ 
                2 6.22\\
                3 3.24\\ 
                4 9.46\\
                };
                \addplot[ybar=0pt, bar width=0.2, fill=orange2, draw=gray1, area legend] table[row sep=crcr] {%
                1 -7.14\\ 
                2 5.08\\
                3 0.54\\ 
                4 1.52\\
                };
                
            \end{axis}
            
            \hfill
            \begin{axis}[%
                width=0.18\textwidth,
                height=0.12\textheight,
                at={(2.9in,0.806in)},
                scale only axis,
                %bar shift auto,
                ybar=1pt,
                log origin=infty,
                xmin=0.509090909090909,
                xmax=4.49090909090909,
                xtick={1,2,3,4},
                ytick={-1,0,1,2,3,4},
                xtick pos=left,
                ytick pos=left,
                xlabel style={font=\color{white!15!black}},
                xlabel={Agent},
                ylabel style={font=\color{white!15!black}, at={(-0.1,0.5)}},
                ylabel={$\Delta_{\%} C_i$},
                %ylabel near ticks,
                label style={font=\footnotesize},
                every tick label/.append style={font=\footnotesize},
                axis background/.style={fill=white},
                legend columns=3,
                legend style={at={(-1.23,1.1)}, anchor=south west, legend cell align=left, align=left, draw=black, font=\footnotesize},
                every axis legend/.append style={column sep=0.1em},
                ]
                
                \addplot[black] table[x index = 0, y index = 1, forget plot] {
                0.5 0
                4.49 0
                };
                
                \draw[ultra thin] (axis cs:\pgfkeysvalueof{/pgfplots/xmin},0) -- (axis cs:\pgfkeysvalueof{/pgfplots/xmax},0);
                
                \addplot[ybar=3pt, bar width=0.2, fill=blue2, draw=gray1, area legend] table[row sep=crcr] {%
                1 -0.8\\ 
                2 0.35\\
                3 2.10\\ 
                4 2.34\\
                };
                \addplot[ybar=3pt, bar width=0.2, fill=new_blue1, draw=gray1, area legend] table[row sep=crcr] {%
                1 -0.8\\ 
                2 0.35\\
                3 1.08\\ 
                4 3.50\\
                };
                \addplot[ybar=3pt, bar width=0.2, fill=orange2, draw=gray1, area legend] table[row sep=crcr] {%
                1 -0.3\\ 
                2 0.28\\
                3 0.18\\ 
                4 0.56\\
                };
                \legend{$\chi=1 \enspace$,$\chi=0$ w/out $K \enspace$,$\chi=0$ w/ $K$}
                
            \end{axis}
    \end{tikzpicture}
    \vspace{-1em}
    \caption{Increase of payments in absolute (left) and percentage (right) values.}
    \label{fig:4_impact}
    \vspace{-1em}
\end{figure}

We hereby propose a small test case to provide in-depth insight on the impacts of different loss allocation policies. We consider a 5 bus system, including both a meshed (TSO) and a radial (DSO) network, with only one generator and three loads, as pictured in \figurename~\ref{fig:4_small_case}. On the one hand, a consumption unit (agent 2) belongs to the TSO grid and is connected to the generator by long lines. On the other hand, agents 3 and 4 are connected at DSO level, respectively close to and far from the feeder. 

We simulate this illustrative example with 3 different loss allocation policies: socialization ($\chi=1$), individual ($\chi=0$ w/out $K$) and individual scaled by capacities ($\chi=0$ w/ $K$). We then investigate the impact of different policies in terms of payment increase in absolute ($\Delta C_i$) and percentage ($\Delta_{\%} C_i$) values,
\begin{equation} \label{eq:delta_payments}
    \Delta C_i = C_i - C_i^0
\end{equation}
\begin{equation} \label{eq:delta_payments_percent}
    \Delta_{\%} C_i = \frac{C_i - C_i^0}{|C_i^0|}
\end{equation}
with $C_i = \sum_{j\in\omega_i} t_{ij} \, (\tau_{ij}^{\textsc{t}} + \tau_{ij}^{\textsc{z}})$ the payments of agent $i$. We use the superscript $0$ to identify the reference case, where losses are not considered while clearing the market.

\figurename~\ref{fig:4_impact} shows that the socialization policy achieves an equal split of losses across market participants. One should note that the revenues of the generator (agent 1) increase, since, in this small case, it is the only generator, hence the marginal one providing the extra power to compensate for the losses that are not considered in the reference case. However, when looking at the results for the same policy in percentage terms, it is clear that the smaller agents (3 and 4) are strongly discriminated. On the other hand, the individual policy allocates larger payments to agents 2 and 4, corresponding to the units further away from the generator. In percentage terms, however, agent 3 is still discriminated over agent 2 because of her size. When including the capacity proportionality to the individual allocation policy, the payments become more uniform across agents in percentage terms. As for the actual payments, they also decrease in magnitude as the marginal generator, having a larger capacity than the others, gets more losses allocated reducing the need for the loads to procure power.
% \vspace{-0.5em}
\section{Numerical Results}\label{sec:4}
This section evaluates the impact of grid constraints and the performance of the different loss allocation policies on a test case, including a TSO network with several connected DSOs. First, the test case is described. Then, the numerical results are investigated to compare market equilibria with and without grid constraints, followed by an analysis on the impact of different loss allocation policies on market outcomes and prices.

\subsection{Test Case Description}

\input{Plots/plot_lineloading.tex}

We simulate our market formulation on a modified RTS `96 IEEE test system presented in \cite{Tosatto2019AAnalysis}. Not being the focus of this paper, the overall synchronicity of the test case is guaranteed by replacing the HVDC line between nodes 123 and 323 by an AC line of the same capacity. To match the formulation of this paper, wind farms are considered as prosumers such that $\underline{P}_i$ and $\overline{P}_i$ are equal to the power production. Focusing on the analysis of the proposed loss allocation policies, the results presented in this paper are obtained for a single time step, September 2$^{\text{nd}}$ of the original dataset. 

\footnotetext[1]{Line capacities of the distribution systems are based on \url{https://portfolio.du.edu/downloadItem/358246}}
\footnotetext[2]{\url{https://gitlab.com/fmoret/p2p_loss}}

\begin{figure}[!b]
	\centering
	\includegraphics[trim={4cm 8cm 1cm 7.5cm},clip,width=0.49\textwidth]{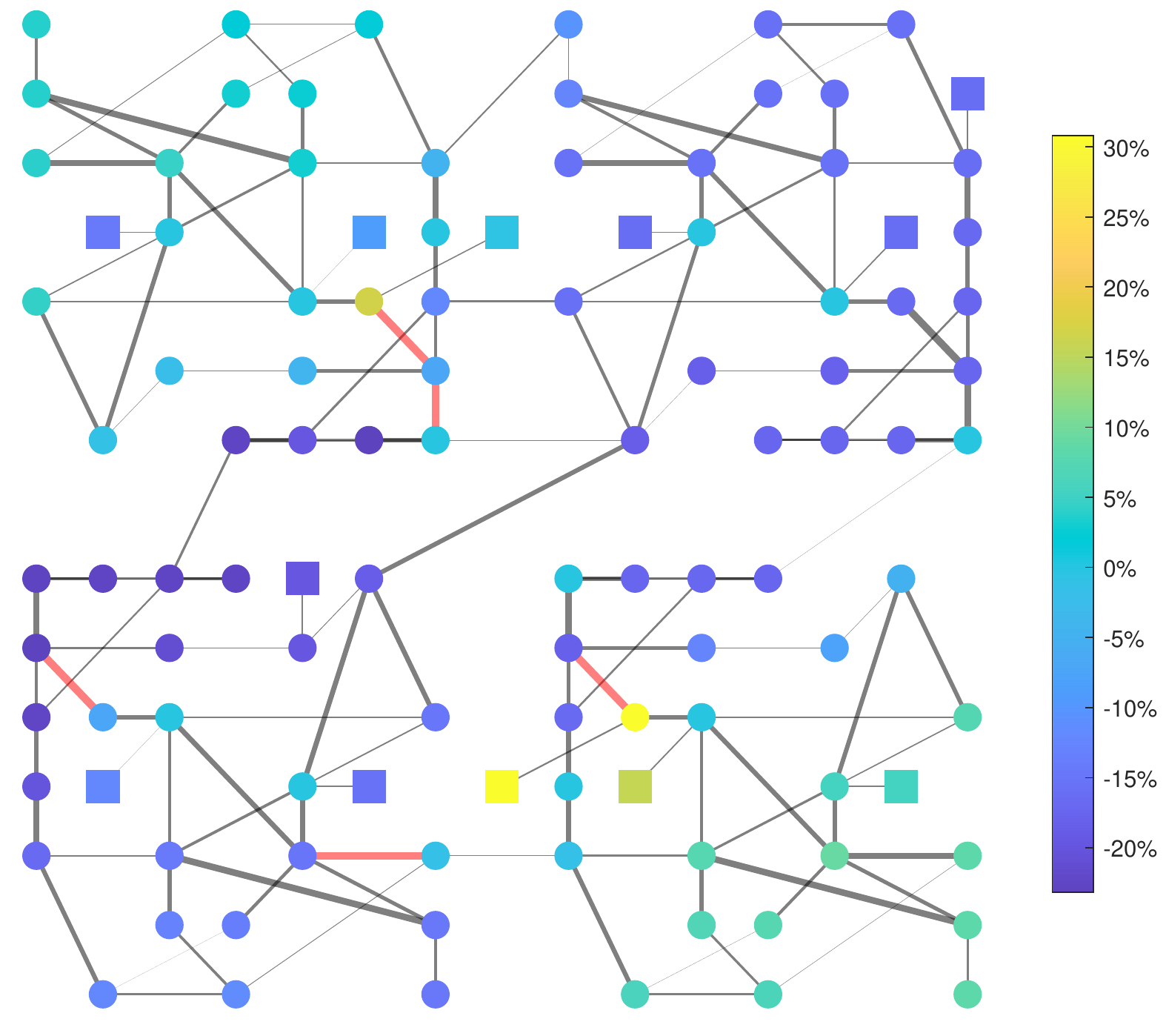}
	\caption{Price percentage increase on the TSO network when including grid constraints.}
	\label{fig:Perceived_Prices}
\end{figure}

Moreover, 12 IEEE 33-bus distribution network\footnotemark[1] connected to transmission nodes 111, 112, 120, 211, 214, 303, 311, 312, 411, 412 and 414 to complete the test case. For the sake of clarity of the results, agents of each distribution networks are gathered in energy communities. However it is worth mentioning that such market structure is only one of the possible ones, as thoroughly discussed in \cite{Baroche2019ProsumerFormulation}.
DSO agents are considered as prosumers with lower bound of their power asset $\bm{\underline{P}}$ corresponding to the power demand of the IEEE 33-bus distribution network (with a change of the sign convention), while their upper bounds $\bm{\overline{P}}$ are computed by randomly sampling the difference $\bm{\overline{P}}-\bm{\underline{P}}$ from a uniform distribution $\mathcal{U}(0,1)$. Similarly, cost curves are sampled uniformly between 10 and 50 \euro/MWh. The final test case as well as the code used to run the simulations are available on GitLab\footnotemark[2]. 
Not being the focus of this paper, simulation results are obtained by solving \eqref{eq:2_6} as a centralized optimization problem. We refer to \cite{Moret2018NegotiationProperties} for further analyses on the computation and communication properties of P2P electricity market solved in a decentralized manner.

\subsection{Impact of Grid Constraints} \label{sec:grid_outcomes}

\begin{figure}[!t]
    \begin{tikzpicture}
        \begin{axis}[%
                width=0.40\textwidth,
                height=0.13\textheight,
                at={(1.20in,0.806in)},
                scale only axis,
                xmin=0,
                xmax=80,
                xlabel style={font=\color{white!15!black}},
                xlabel={Weighted electrical distance ($\delta$)},
                ymin=0,
                ymax=18,
                ylabel style={font=\color{white!15!black}},
                ylabel={Number of agents},
                ytick={0,5,10,15,20,25},
                axis background/.style={fill=white},
                legend columns=2,
                legend style={at={(0.28,1.05)}, anchor=south west, legend cell align=left, align=left, draw=black, font=\footnotesize},
                every axis legend/.append style={column sep=0.1em}
                ]
                
                \addplot[ybar interval, fill=new_orange1, fill opacity=0.6, draw=gray1, area legend] table[row sep=crcr] {%
                x	y\\
                1	1\\
                3	0\\
                5	0\\
                7	0\\
                9	1\\
                11	1\\
                13	2\\
                15	3\\
                17	2\\
                19	6\\
                21	9\\
                23	5\\
                25	3\\
                27	7\\
                29	1\\
                31	0\\
                33	1\\
                35	0\\
                37	3\\
                39	6\\
                41	10\\
                43	4\\
                45	9\\
                47	16\\
                49	7\\
                51	5\\
                53	5\\
                55	4\\
                57	4\\
                59	3\\
                61	5\\
                63	3\\
                65	4\\
                67	1\\
                69	5\\
                71	1\\
                73	1\\
                75	6\\
                77	0\\
                79	0\\
                };
                \addlegendentry{No grid}
                
                \addplot[ybar interval, fill=new_green1, fill opacity=0.6, draw=gray1, area legend] table[row sep=crcr] {%
                x	y\\
                1	17\\
                3	7\\
                5	9\\
                7	5\\
                9	14\\
                11	10\\
                13	13\\
                15	13\\
                17	5\\
                19	6\\
                21	8\\
                23	5\\
                25	9\\
                27	5\\
                29	7\\
                31	1\\
                33	1\\
                35	15\\
                37	1\\
                39	3\\
                41	1\\
                43	2\\
                45	0\\
                47	0\\
                49	1\\
                51	1\\
                53	0\\
                55	1\\
                57	1\\
                59	0\\
                61	0\\
                63	0\\
                65	0\\
                67	0\\
                69	0\\
                71	0\\
                73	0\\
                75	0\\
                77	0\\
                79	0\\
                };
                \addlegendentry{Grid}

            \end{axis}
    \end{tikzpicture}
    \vspace{-1em}
    \caption{Weighted average electric distance of trades per TSO prosumer.}
    \label{fig:Trading_Distances}
    \vspace{-0.2em}
\end{figure}
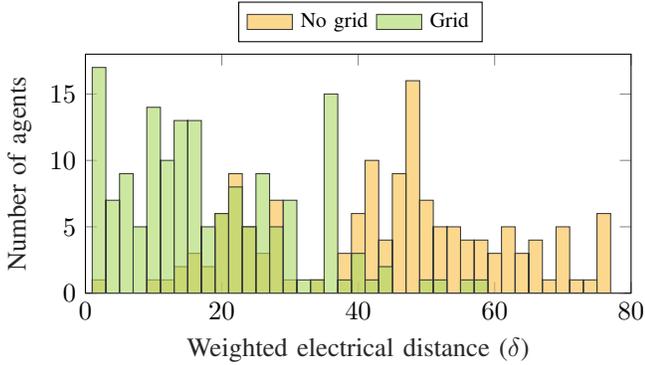
\input{Plots/plot_nocap.tex}

We first investigate the impact of power flow constraints on the feasibility of market outcomes, by simulating the market clearing \eqref{eq:2_6} with and without including grid constraints. This analysis highlights the advantage of the endogenous P2P market proposed here compared to an exogenous approach such as in \cite{Baroche2019ExogenousMarketsb}. Power losses and their allocation policies are not yet considered, to only analyze the effect of the network on the market equilibrium. The results of our simulations verify that all power grid variables respect their limits, once included in the model. For instance, \figurename~\ref{fig:lines_loading} shows the loading of transmission and distribution lines for a market clearing with and without grid constraints. The results clearly show violation of line capacities for both system operators, if grid constraints are not included in the negotiation process, in particular at DSO level where reactive flows contribute to the overload of lines. Optimizing energy trades while accounting for grid feasibility become fundamental in view of an even lager penetration of DERs, causing considerable reversed (from DSO to TSO) and reactive power flows. By including grid constraints in the market clearing, system operators will reduce the need for grid reinforcement and reactive compensation measures.

From \eqref{eq:KKT}, we extracted how grid constraints modify trade prices to account for power flow feasibility, e.g., line and voltage limits. Similarly to a market based on locational marginal pricing, one line congestion impacts the price of each bus. Note that agents connected at the same bus trade at the same price if no economic preference is included. \figurename~\ref{fig:Perceived_Prices} displays the percentage increase, compared to the market clearing without grid constraints, of the price of agents at each TSO bus. For clarity of representation, DSOs (represented with a squared marker) are grouped at their connection node to the TSO grid. The width of the connecting edges is plotted proportionally to the line flows, with congested lines highlighted in red. The average price per agent decreases from 15.6 \euro/MWh (the uniform price of the market clearing without grid constraints) to 14.5 \euro/MWh. However, larger prices, with a peak of 20.4 \euro/MWh, appear as a consequence of congested lines. Visually, it appears that prices are distributed in high and low price zones: in fact, including grid constraints in the market clearing provides incentives to locally match generation and consumption. 

\input{Plots/plot_lossTSODSO.tex}

To this extent, \figurename~\ref{fig:Trading_Distances} shows how the average trading distance of TSO prosumers is largely reduced both in terms of its mean and standard deviation. DSO agents are not considered, since, in this specific market architecture, they only trade with the respective community manager located at the feeder, hence their distance is not impacted. We calculate the average trading distance $\delta_i$ of each agent $i$ as the PTDF-based electric distance between two trading agents weighted by the energy traded. It follows
\begin{equation} \label{eq:el_dist}
    \delta_i = \frac{\sum_{j\in\Omega_i} \sum_{l\in\mathcal{L}} |\text{PTDF}_{li} - \text{PTDF}_{lj}| \, r_l \, t_{ij}}{\sum_{j\in\Omega_i} t_{ij}}
\end{equation}
with $r_l$ the resistance of each line $l$. As investigated in \cite{Moret2018NegotiationProperties}, the number of trades of each market participant increases the complexity of the market clearing problem, especially if solved by means of distributed optimization techniques. The natural effect of the power grid to incentivize local trades presented here could be employed to increase the sparsity of the communication matrix among agents.

\subsection{Impact of Loss Allocation Policies} \label{sec:loss_outcomes}

\input{Plots/plot_cap.tex}

After investigating the impact of grid constraints on market outcomes, we focus our analysis on the different loss allocation policies proposed. We first compare the effects of a socialization \eqref{eq:soc_pol} and an individual \eqref{eq:ind_pol} policy. \figurename~\ref{fig:deltaP_dist} displays the percentage increase of payments as a function of the weighted electrical distance, defined in \eqref{eq:el_dist}, of both TSO and DSO agents. We use the socialization policy as a reference, where each trade is allocated an equal share of the respective SO losses, and plot the cost increase when an individual policy is employed. It is clear that market participants that engage in distant trades are penalized by larger losses, and consequently payments. This trend is stronger at DSO rather than at TSO level, since agents at the end of a radial grid are naturally further away from the feeder, while it is attenuated in a meshed grid where power flows through several parallel paths.

This can be clearly seen from \figurename~\ref{fig:loss_dist}, where losses and electrical distance are displayed for a sample of the simulated agents. While a socialization policy leads to equal losses across all SO agents, an individual policy induces large disparity at DSO level as function of their electrical distance from the feeder. On the other hand, market participants at TSO level are indeed impacted by their trading electrical distance, but the variation is limited thanks to the meshed layout of the power grid. Therefore, we argue that a socialization policy is preferred at DSO level to avoid discriminating agents for their unintentional geographical location, while at TSO level an individual policy can still be used, charging more for trades with higher grid usage, but without creating unfairly large payment differences.

As pointed out in Section \ref{sec:3}, an ideal individual policy should use the actual amount of trades to quantify losses, but this would introduce bilinear terms, mining the convexity of the overall problem. We propose to consider the capacity of market participants as a proxy of their trades. \figurename~\ref{fig:deltaP_cap} shows the impact of considering agent capacities in an individual allocation policy, as in \eqref{eq:cap_pol}, by displaying the percentage payments increase (with the socialization policy as a reference) as function of the overall energy traded by each agent. The clear linear trend obtained when including agent capacities in the allocation policy is evidence that we can approximate the objective of allocating more losses to agents that trade more energy, while keeping the optimization problem convex.

% \vspace{-0.5em}

%% CONCLUSION %%%%%%%%%%%%%%%%%%%%%%%%%%%%%%%%%%%%%%%%%%%%%%%%%%%%%%%%%%%%%%%%%%%%%%%%%%%%%%%%%%%%%%%%%%%%%%%%%%%%%%
\section{Conclusion}\label{sec:5}
In view of decentralized electricity markets, the role of system operators is fundamental to guarantee feasibility in power system operation. In this work, we tackled the coordination of system operators with agents participating in peer-to-peer electricity markets and with other system operators, both at transmission and distribution level. We showed that by including system operators in the negotiation process grid constraints are respected and trade prices are affected by network charges. These additional tariffs are dynamically calculated as an outcome of the negotiation mechanism, reflecting the actual grid usage of each trade and providing incentives to local match of generation and consumption. We then included power losses as an additional market product that system operators trade for compensating transmission and distribution costs. We finally argued that loss allocation policies are needed to guarantee fairness among market participants, that could be discriminated for their geographical location otherwise.

In order to focus the analysis and contribution solely on market design, simplified convex approximations of the power flow dynamics were employed. Since the optimization problems of system operators are somewhat separated from those of other market agents (in the sense that, market mechanisms often simplify power system operation), the proposed approach and related market mechanisms will not be affected by more complex power flow modelling as long as it is convex. Further research should tackle such improvement, by enhancing the accuracy of the power flow and loss modelling. Furthermore, this work laid out the foundation for market-based integration of system operators in decentralized electricity markets, including market products for allocating power losses to market participants. More advanced loss allocation policies could be developed as future work, but they could still be implemented by means of the same market products proposed in this work.

Including system operator in the negotiation mechanism, based on decomposition techniques for large optimization problems, drastically increases the complexity of the algorithm. Even if it is still possible to find the optimal consensus among all these market actors, the time needed to attain it might increase to the point of making such mechanism not suitable for real life implementations. Therefore, this work should be considered as an ideal benchmark for the design of decentralized electricity markets. For instance, it could be used for exploring new opportunities for a decentralized management of the power systems, as well as for better defining \textit{ex-ante} network charges that approximate the results proposed in this work. Finally, further research should aim for designing negotiation processes with guaranteed trade-off between optimality of market equilibria and algorithmic complexity. Computationally lighter mechanisms may lead to market inefficiency, but this might be the price to pay for real-life implementation of such complex market mechanisms.

%% BIBLIOGRAPHY %%%%%%%%%%%%%%%%%%%%%%%%%%%%%%%%%%%%%%%%%%%%%%%%%%%%%%%%%%%%%%%%%%%%%%%%%%%%%%%%%%%%%%%%%%%%%%%%%%%%
\bibliographystyle{myIEEEtran.bst}

%% BIOGRAPHY %%%%%%%%%%%%%%%%%%%%%%%%%%%%%%%%%%%%%%%%%%%%%%%%%%%%%%%%%%%%%%%%%%%%%%%%%%%%%%%%%%%%%%%%%%%%%%%%%%%%
\begin{IEEEbiography}[{\includegraphics[width=1in,height=1.25in,clip,keepaspectratio]{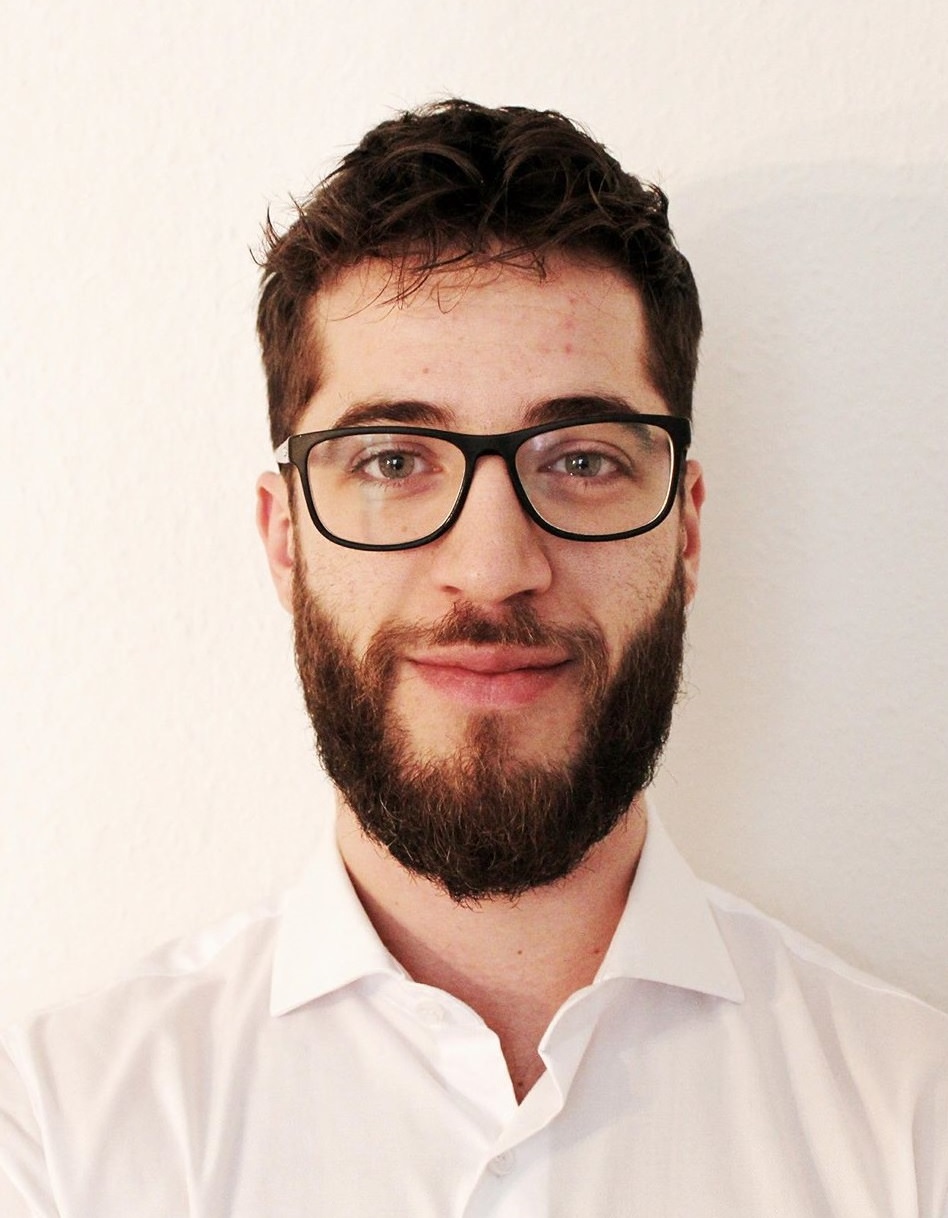}}]{Fabio Moret} received the M.Sc. degree in electrical engineering from the University of Padova, Padova, Italy, and the M.Sc. degree in sustainable energy from the Technical University of Denmark, Lyngby, Denmark, where he received a Ph.D. degree at the Centre for Electric Power and Energy, Department of Electrical Engineering. He is currently a data scientist at A.P. Moller – Maersk in the Uptake Management department. His research interests include consumer-centric electricity markets, distributed negotiation mechanisms, decomposition techniques for optimization problems, and applied game theory.
\end{IEEEbiography}
\vspace{-1em}
\begin{IEEEbiography}[{\includegraphics[width=1in,height=1.25in,clip]{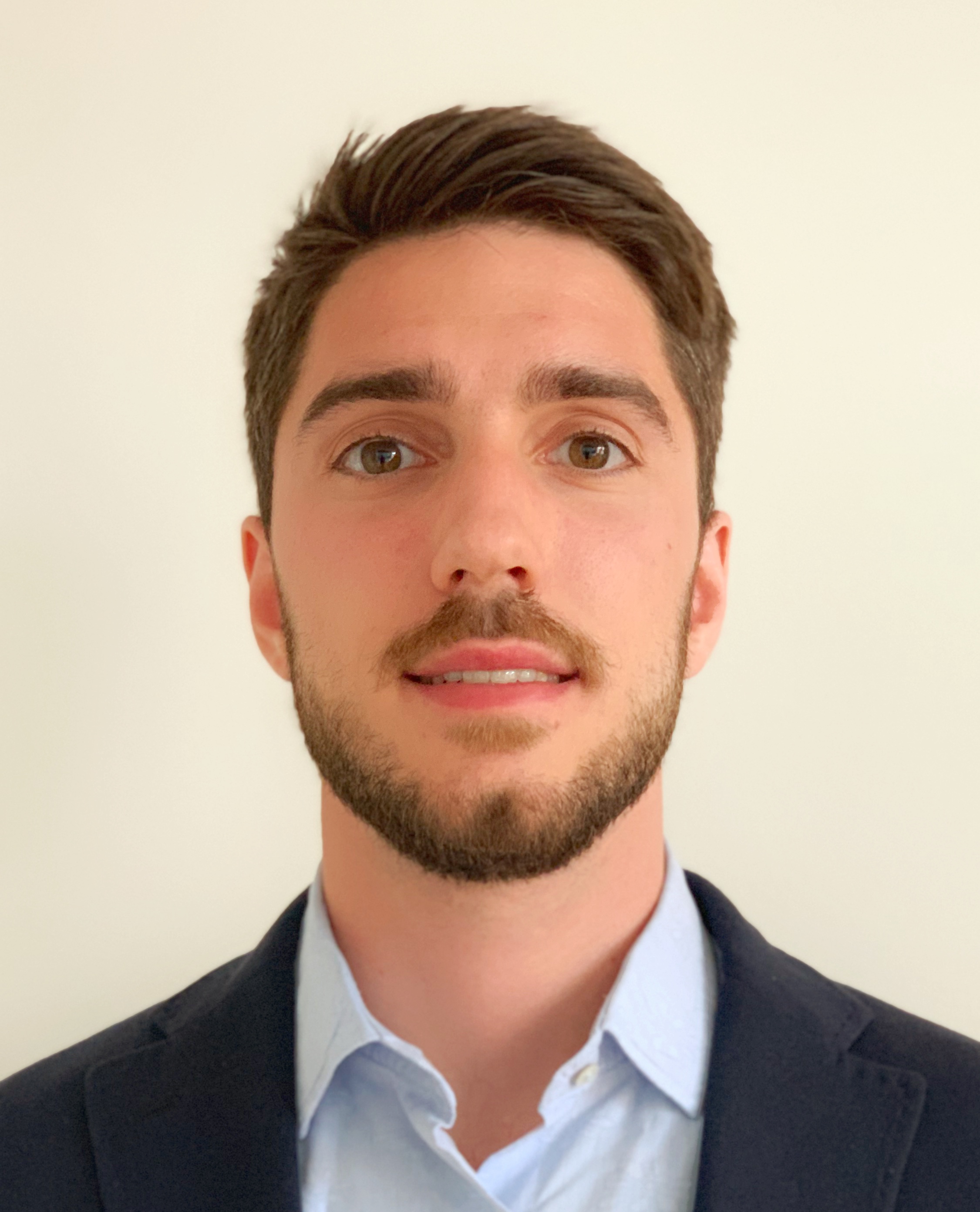}}]{Andrea Tosatto}
(S'18) received the M.Sc. degree in electrical engineering from the Royal Institute of Technology, Stockholm, Sweden, and a second M.Sc. degree in electrical engineering from the Institute Polytechnique de Grenoble, Grenoble, France. He is working toward the Ph.D. degree at Centre for Electric Power and Energy, Department of Electrical Engineering, Technical University of Denmark. His research interests include convex optimization in power systems, applied game theory, and market integration of multiarea AC/HVDC grids.
\end{IEEEbiography}
\vspace{-1em}
\begin{IEEEbiography}[{\includegraphics[width=1in,height=1.25in,clip,keepaspectratio]{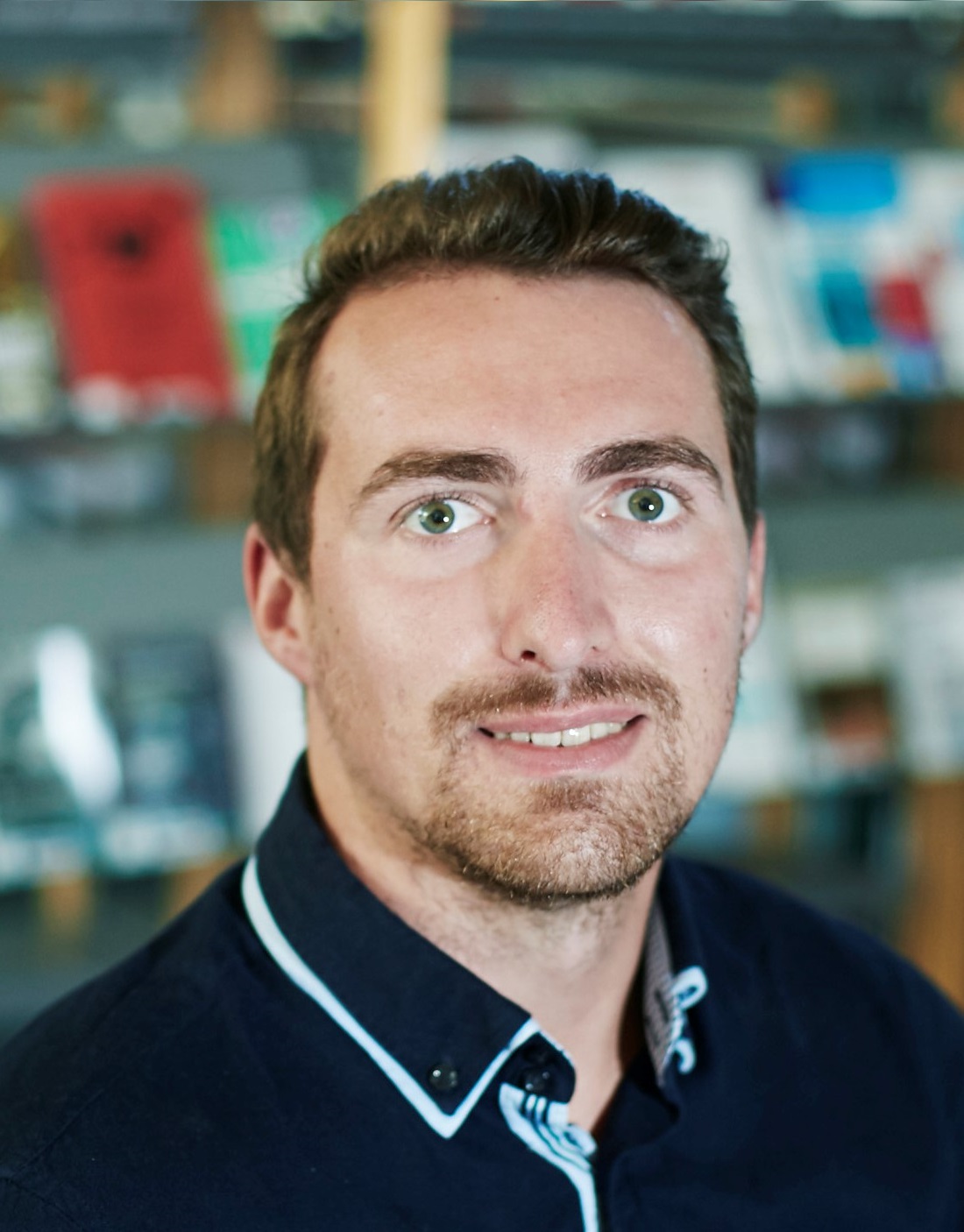}}]{Thomas Baroche} received the M.Sc. degree in electrical engineering from the Ecole Normale Supérieure de Rennes, Rennes, France and the M.Sc. degree in renewable energy, sustainability and technology from the Ecole Polytechnique, Palaiseau, France. He is working toward the Ph.D. degree in electrical engineering at Ecole Normale Supérieure de Rennes, Rennes, France. His research interests include decentralized peer-to-peer electricity markets and grid constraints integration. 
\end{IEEEbiography}
\vspace{-1em}
\begin{IEEEbiography}[{\includegraphics[width=1in,height=1.25in,clip]{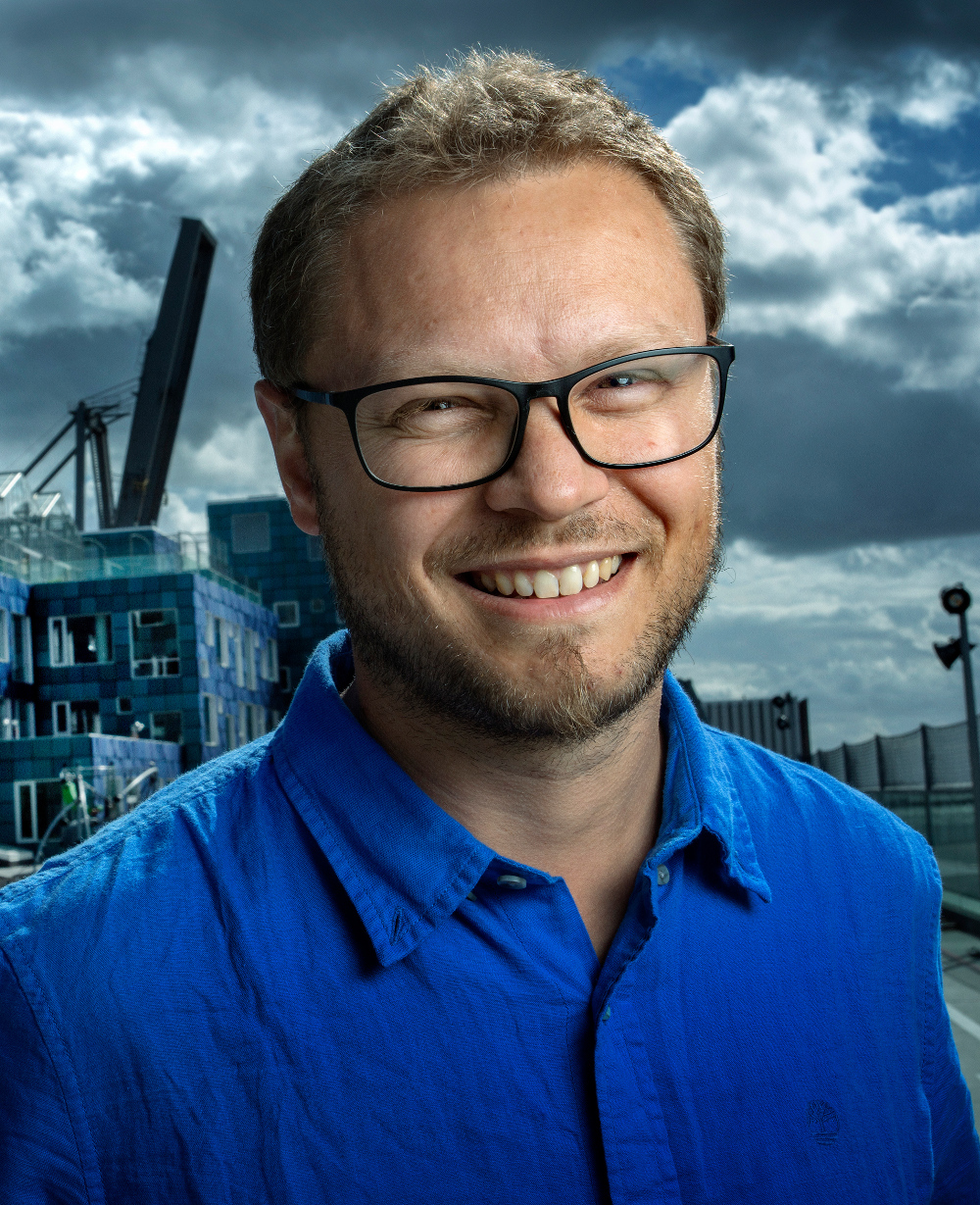}}]{Pierre Pinson}
(SM’13, F’20) received the M.Sc. degree in applied mathematics from the National Institute for Applied Sciences, Toulouse, France, and the Ph.D. degree in energetics from Ecole des Mines de Paris, Paris, France. He is a Professor with the Centre for Electric Power and Energy, Department of Electrical Engineering, Technical University of Denmark, Lyngby, Denmark, also heading a group focusing on energy analytics and markets. His research interests include operations research and management science, with application within power and energy systems among others. He is the Editor-in-Chief for the International Journal of Forecasting.
\end{IEEEbiography}

\vfill

\end{document}